\documentclass[12pt]{article}

\usepackage[a4paper, total={6.4in, 8.8in}]{geometry}
\usepackage{amsmath,amssymb,latexsym,amsfonts}
\usepackage{stmaryrd}
\usepackage{shapepar,bm}
\usepackage[dvips,arrow,matrix,ps,color,line]{xy}
\usepackage{theorem,amsfonts,amssymb,amscd}
\usepackage{geometry}
\usepackage{makeidx}
\usepackage{eulervm}

\usepackage{graphicx}
\usepackage{graphics}

\usepackage{epstopdf}

\usepackage{amssymb,graphics,hyperref}

\hypersetup{
    colorlinks = true,
    linkcolor = blue,
    anchorcolor = red,
   citecolor = blue,
    filecolor = red,
    pagecolor = red,
    urlcolor = blue}

\newcommand{\mathsym}[1]{{}}

\usepackage[usenames]{color}
\definecolor{MyLightMagenta}{cmyk}{0.1,0.8,0,0.1}
\definecolor{MyDarkBlue}{rgb}{0.1,0,0.3}

\makeindex

\def\bft{\mathbf{t}}

\def\NN{\mathbb N}

\def\ubX{\mathbf{X}}

\def\bfz{{\mathbf z}}

\def\ovsig{\overline{\sigma}}

\def\DJKM{Date, Jimbo, Kashiwara and Miwa}

\def\bfX{{\mathbf X}}

\def\bfu{{\mathbf u}}

\def\d{\partial}

\def\Vcal{\mathcal V}

\def\bfx{{\mathbf x}}

\def\bft{{\mathbf t}}
\def\bfs{{\mathbf s}}

\def\wX{[{\mathbf X}]}

\def\ZZ{\mathbb Z}

\def\CC{\mathbb C}

\def\Ecal{{\mathcal E}}

\def\QQ{\mathbb Q}
\def\PP{\mathbb P}

\def\cocoa{{\hbox{\rm C\kern-.13em o\kern-.07em C\kern-.13em o\kern-.15em A}}}

\def\Dcal{\mathcal D}

\def\Fcal{{\mathcal F}}

\def\bmd{{\bm\d}}

\def\bft{{\bf t}}

\def\bfu{{\bf u}}

\def\End{\mathrm{End}}

\def\blamb{{\bm \lambda}}
\def\bmu{{\bm \mu}}

\def\Pcal{{\mathcal P}}

\def\w2M{\bigwedge^2M}

\def\w{\wedge }
\def\bw{\bigwedge }
\def\bwV{{\bigwedge\hskip-3.5pt V}}

\protect
\protect
\protect
\protect

\protect
\protect

\def\sra{\rightarrow}

\def\proof{\noindent{\bf Proof.}\,\,}
\def\qed{{\hfill\vrule height4pt width4pt depth0pt}\medskip}
\def\be{\begin{equation}}
\def\ee{\end{equation}}
\def\bclm{\begin{claim}}
\def\eclm{\end{claim}}
\def\beqn{\begin{eqnarray}}
\def\eeqn{\end{eqnarray}}
\def\beqn*{\begin{eqnarray*}}
\def\eeqn*{\end{eqnarray*}}

\theoremstyle{change}
\theorembodyfont{\rmfamily}

\newtheorem{claim}{}[section]

\global
\global

\def\no@breaks#1{{\def\\{ \ignorespaces}#1}}    

  \makeatletter
\def\cleardoublepage{\clearpage\if@twoside \ifodd\c@page\else
\hbox{} \thispagestyle{empty}
\newpage
\if@twocolumn\hbox{}\newpage\fi\fi\fi} \makeatother

\usepackage{eso-pic,graphicx}
\makeatletter
\newcommand\BackgroundPicture[2]{%
  \setlength{\unitlength}{1pt}%
  default \put(0,\strip@pt\paperheight){%
  \parbox[t][\paperheight]{\paperwidth}{%
    \vfill
     \centering \includegraphics[angle=#2, width=15cm, height=15cm,  bb=0 0 150 150]{#1}
    \vfill
}}} %
\makeatother

\newcommand{\llb}{\llbracket}
\newcommand{\rrb}{\rrbracket}
\def\bmtx{\begin{matrix}}
\def\emtx{\end{matrix}}
\makeindex

\begin{document}

\date{}
\title{On the Vertex Operator Representation of Lie Algebras of Matrices}

\author{Ommolbanin Behzad, Andr\'e Contiero, \\ 
David Martins
\thanks{The second author is partially supported by Funda\c c\~ao de Amparo \`a Pesquisa do Estado
de Minas Gerais (FAPEMIG) grant no. APQ-00798-18. The third author was financed in part by the 
Coordena\c c\~ao de Aperfei\c coamento de Pessoal de N\'ivel Superior - Brasil (CAPES) - Finance Code 001.
\smallskip 
\newline  ${}$
\,\,\,\,\,\,\,{\em Keywords and Phrases:} Hasse-Schmidt Derivations, Vertex Operators  on Exterior Algebras,  Representation of Lie Algebras of Matrices,
Bosonic and Fermionic Representations by Date-Jimbo-Kashiwara-Miwa, Symmetric Functions.
\smallskip
\newline ${}$ \,\,\,\,\,\, {\bf 2020 MSC:} 17B10, 17B69, 14M15, 15A75,  05E05.}}

\maketitle

\vspace{-28pt}
\begin{abstract}

\noindent
The polynomial ring $B_r:=\QQ[e_1,\ldots,e_r]$ in $r$ indeterminates is a 
representation of the Lie algebra of all the endomorphism of $
\QQ[X]$ vanishing at powers $X^j$ for  all but finitely many $j$. We determine a $B_r$-valued formal power series in $r+2$ indeterminates which encode the images of all the basis elements of $B_r$ 
under the action of the generating function of elementary endomorphisms of $\QQ[X]$, which we call the {\em structural series of the representation.} The obtained expression  implies (and improves) a formula by Gatto \& Salehyan,  which only computes, for one chosen basis element, the  generating 
function of its images. For sake of completeness we construct in the last section  the  $B=B_\infty$-valued structural formal power series which consists in the evaluation of the  vertex 
operator describing the bosonic representation of $gl_\infty(\QQ)$ against the generating function of the standard Schur basis of $B$. This provide an alternative  description of the bosonic representation of $gl_\infty$ due to Date, Jimbo, Kashiwara and Miwa which does not involve explicitly exponential of differential operators. 

\end{abstract}


\section{Introduction}
\claim{}
This paper is concerned with the following general and rather elementary 
fact. There is a natural way to multiply polynomials by matrices of infinite 
size (with all zero entries but finitely many) which  is compatible with the  
matrix commutator, i.e. $M(Np)-N(Mp)=[M,N]p$ for all matrix pair $(M,N)$ and 
for all polynomials  $p$.  This observation is basically due to \DJKM \, 
(DJKM) \cite{DJKM01} who, more than that, determine a vertex operator 
representation of the generating function of the elementary matrices acting 
on the ring $B$ (the bosonic Fock space) of polynomials in infinitely many 
indeterminates  (see also \cite[Section 5.1]{KacRaRoz} for an elementary 
account).

Our main contribution consists in the study of the {\em finite type} version 
of the DJKM description of the bosonic representation of matrices, improving 
the output of the approach taken in \cite{gln}. To be more precise  we need 
to introduce two main actors. 

The former is the vector space $V:=\QQ[X]$  of polynomials in the 
indeterminate $X$, with basis $
(X^i)_{i\geq 0}$, to which one attaches the Lie algebra $gl(V):=
\bigoplus_{i,j\geq 0}\QQ\cdot X^i\otimes\d^j$, where  $\d^j$ denotes the 
unique linear form such
that $\d^i(X^j)=\delta^{ij}$. 

The latter is the polynomial ring $B_r:=\QQ[e_1,\ldots,e_r]$,  
 in the $r$ indeterminates $(e_1,\ldots,e_r)$. Denote by $H_r$ the sequence $
 (h_n)_{n\in\ZZ}
$ in $B_r$, where $h_n:=\det(e_{i-j+1})_{1\leq i,j\leq n}$, setting by 
convention, $e_0=1$ and $e_j=0$ for $j<0$. Let $\Pcal_r$ be the set of all 
partitions of length at most $r$. It turns out that the set of all $\Delta_
\blamb(H_r)=\det(h_{\lambda_j-j+i})$ provide a $\QQ$-basis of $B_r$  
parametrized by $\Pcal_r$, i.e. $B_r:=\bigoplus_{\blamb\in\Pcal_r}\QQ\cdot 
\Delta_
\blamb(H_r)$. The main point is that $B_r$ can be made into an irreducible 
representation of $gl(V)$  (in particular: polynomials in $B_r$ can be 
multiplied by matrices  compatibly with the Lie algebra structure of $gl(V)
$), by pulling back the natural one on $\bw^rV$  through the vector space 
isomorphism $B_r\sra \bw^rV$ as in \cite[Main Theorem]{LakTh01}, the {\em 
finite type boson-fermion correspondence.}

 Our Theorem \ref{thm:mnthm1} determines an expression  for the formal power 
 series $\Ecal_r(z,w^{-1},\bft_r)\in B_r\llb z,w^{-1},\bft_r\rrb$ resulting 
 from the evaluation of the generating function
$$
\Ecal(z,w^{-1})=\sum_{i,j\geq 0}X^i\otimes \d^j\cdot z^iw^{-j}
$$ 
of the  elementary endomorphisms $X^i\otimes \d^j\in \End_\QQ(V)$ against a 
suitable generating function  of the basis $\big(\Delta_
\blamb(H_r)\big)_{\blamb\in\Pcal_r}$ of $B_r$.

For sake of a dutiful comparison,  we apply the same procedure 
to the genuine DJKM- representation of $gl_
\infty(\QQ)$ in Section \ref{sec:sec5}, by  merely evaluating  the bosonic 
vertex operator against the generating function of 
 the natural Schur basis of $B$, which can be determined via Cauchy-type 
 formulas. One so  obtains  a description of the $gl_\infty(\QQ)$--module 
structure of $B=B_\infty$ (Theorem~\ref{thm:mnthm2}) which is equivalent to 
the DJKM one, but with no explicit occurrence of exponentials of differential 
operators. The rest of the introduction will be 
devoted to state more precisely our results and to say a few words about the 
history of the subject and motivations.

\claim{\bf Precise Statement of the Main Results.} The most natural candidate 
to be a generating series for the 
basis $(\Delta_\blamb(H_r))$ is 
$$
\Delta(H_r;\bft_r):=\sum_{\blamb\in\Pcal_r}\Delta_\blamb(H_r)s_\blamb(\bft_r)
$$
where $\bft_r:=(t_1,\ldots,t_r)$ is an $r$-tuple of indeterminates and 
$s_\blamb(\bft_r)$ are the usual Schur symmetric polynomials as defined, 
e.g., in \cite[Formula 3.1]{MacDonald}. 
The finite type version of the {\em boson fermion correspondence} says that 
the $r$-th exterior power $\bw^rV$ is a free $B_r$--module of rank $1$  
generated by $X^r(0):=X^{r-1}\w\cdots\w X^0$ such that 
$$
X^r(\blamb):=X^{r-1+\lambda_1}\w\cdots\w X^{\lambda_r}=\Delta_\blamb(H_r)
\cdot X^r(0).
$$
One then defines a $\star$--product, to make $B_r$ into a $gl(V)$-module, as 
follows:
$$
[(X^i\otimes \d^j)\star \Delta_\blamb(H_r)]X^r(0):=X^i\w (\d^j\lrcorner 
X^r(\blamb))
$$ 
where $\d^j\lrcorner :\bw V\sra \bw V$ is the natural contraction operator. 
Define 
\be\label{Erw}
E_r(w)=1-e_1w+\cdots+(-1)^re_rw^r\in B_r[w]
\ee
and 
\be\label{Erwtr}
E_r\left(\bft_r;\displaystyle{1\over w}\right)=1-e_1(\bft_r)w^{-1}+\cdots+
(-1)^re_r(\bft_r)w^{-r}=\prod_{j=1}^r(1-t_jw^{-1})\in \QQ[\bft_r,w^{-1}]
\ee
in such a way that $e_j(\bft_r)$ is the elementary symmetric polynomial of 
degree $j$ in the variable $(t_1,\ldots, t_r)$. The polynomials $e_j(\bft_r)$ 
should not be confused with the indeterminates $e_j$ of the ring $B_r$, which 
could be rather thought of as the elementary symmetric polynomials of the 
universal roots of the generic monic polynomial $X^r-e_1X^{r-1}+\cdots+
(-1)^re_r$. 
Our main result consists in computing the action of the generating function $
\Ecal(z, w^{-1})$ against the generating function $\Delta(H_r;\bft_r)$ of the 
basis of $B_r$.

\noindent
\medskip
{\bf Theorem A} (Theorem \ref{thm:mnthm1}){\bf .} {\em  Let
$$\Ecal_r(z,w^{-1},\bft_r):=\sum_{((i,j),\blamb)\in\NN^2\times\Pcal_r}\left[X^i\otimes \d^j\star \Delta_\blamb(H_r)\right]z^iw^{-j}s_\blamb(\bft_r)
$$
be the generating function of the structure constant of $B_r$ as a representation of $gl(V)$. Then:
$$
\Ecal_r(z,w,\bft_r)=
$$
\smallskip
\begin{equation*}
    {z^{r-1} \over w^{r-1}}\exp\left(\sum_{n\geq 1}{1 \over n}p _n(\bft_r)\left({1 \over w^n}-{1 \over z^n}\right)+x_np_n(z,\bft_r)\right)\left(E_r(w)+(-1)^{r+1}e_rw^rE_r\left(\bft_r,\displaystyle{1\over w}\right)\right).
\end{equation*}
}
The proof of Theorem A  is based on the formalism of Schubert derivations, in the same vein of \cite{BeCoGaVi,BeGa,BeNa}. 

It is then  natural, beside being dutiful, to see how the same  procedure can be applied to rephrase the DJKM-description of the $gl_\infty(\QQ)$-structure of $B=B_\infty$. In this case we identify $B$ with the polynomial ring $\QQ[x_1,x_2,\ldots]$, where the variables $\bfx$ are related to the  $e_i$ through the equality
$$
\exp(\sum_{i\geq 1}x_iz^i)={1\over E_\infty(z)}=(1-e_1z+e_2z^2+\cdots)^{-1}\in B\llb z \rrb.
$$
To be more adherent with the standard notation in the case of $B_\infty$, the Schur basis will be denoted by $S_\blamb(\bfx)=\det(S_{\lambda_j-j+i}(\bfx))$ where $S_k(\bfx)=\det(e_{i-j+1})_{1\leq i,j\leq k}$. The generating function of the basis elements of $B$ is then $\sum_{\blamb\in\Pcal}S_\blamb(\bfx)s_\blamb(\bft)$, where $\bft:=\bft_\infty$. Then our second main result  is:

\medskip
\noindent
{\bf Theorem B} (Theorem~\ref{thm:mnthm2}){\bf .}  {\em
Let 
$$
\Ecal(z,w^{-1}, \bft_r)=\sum_{((i,j,\blamb)\in\ZZ\times \ZZ\times \Pcal}(X^i\otimes \d^j)\star S_\blamb(\bfx))z^iw^{-j}s_\blamb(\bft)\in B \llb z,w, \bft_r,z^{-1}, w^{-1}\rrb 
$$ 
where $\bft:=(t_1,t_2,\ldots)$ is a sequence of infinitely many indeterminates.
Then 
\begin{equation}
\Ecal(z,w^{-1}, \bft_r)=
\exp \left(\sum_{n\geq 1}{1\over n}\left({w^n \over z^n} -{p_n (\bft_r) \over z^n}+{p_n (\bft_r) \over w^n}\right)
+x_n \left(z^n -w^n+p_n (\bft_r)\right)
\right)\label{eq:thmB}
\end{equation}
i.e. the image of $S_\blamb(\bfx)\in B$ through the multiplication by the elementary endomorphism $X^i\otimes \d^j$ is he coefficient of $z^iw^{-j}$ of \eqref{eq:thmB}.
}

\claim{\bf History and motivations.} The bosonic vertex operator representation of the Lie algebra $gl_\infty(\CC)$ was determined by \DJKM \, within the framework of algebraic analysis and mathematical physics related to the KP hierarchy, see e.g. \cite{jimbomiwa}. The KP hierarchy is a system of infinitely many non linear PDEs whose polynomial solutions are parametrized by the points of the orbit of $1\in B\otimes_\QQ\CC$ through a natural action of the group
\begin{center}
$GL_\infty(\CC):=\{$invertible $A\in \End_\CC(\CC[X^{-1},X])\,|\, AX^j=X^j$ for all but finitely many $j\}.$
\end{center}
The $GL_\infty(\CC)$ orbit of $1$ does correspond to the locus of the 
decomposable tensors (the Sato Universal Grassmann Manifold, see 
\cite{SatoUGM} and also \cite[p.~73]{KacRaRoz}) in infinite wedge power, 
roughly speaking $\bw^\infty\CC^\infty$, the Fermionic Fock space of charge 
$0$, to which $B$ is isomorphic via the so called boson-fermion 
correspondence. It turns out that the DJKM representation is the 
linearization of this natural action of $GL_\infty(\CC)$ on $B$. 

In \cite{gln,SDIWP} is recognized that the DJKM description is a natural 
consequence of the well known basic linear algebraic fact, namely that each 
vector space is a representation of its Lie algebra of endomorphisms. More 
generally,  it turns out that the DJKM representation is a particular 
(extremal) case of a more general picture which in \cite{gln} was summarized 
by the slogan ``the cohomology of the Grassmannian is a $gl_n$--module'', 
which entitles the paper. 

This can be quickly explained as follows. Each $\QQ$-vector space $V$ of 
finite dimension $n$ is isomorphic to $£\QQ[X]/(X^n)=H^*(\PP^{n-1},\QQ)$ and 
one already sees that the singular cohomology of the projective space $
\PP^{n-1}$ is a module over the Lie algebra of $\QQ$-valued square $n\times n
$ matrices. The trivial claim that $\QQ[X]$ is a $gl(\QQ[X])$-module 
generalizes to the fact that $B_r$ is a $gl(\QQ[X])$-module because $\bw^r
\QQ[X]$ is naturally a representation of $gl(\QQ[X])$ and because of the 
isomorphism $B_r\sra \bw^r\QQ[X]$. The composition of maps $B_r\sra \bw^r
\QQ[X]\sra \bw^r\QQ[X]/(X^n)$ factorizes through a ring $B_{r,n}$ which turns 
out to be the  the singular cohomology ring of the complex Grassmannian 
$G(r,n)$. The case $n=\infty$ corresponds to the situation coped with in our 
Theorem \ref{thm:mnthm1}. 

One main point is that we do not know any direct way to infer our Theorem~
\ref{thm:mnthm1} from the DJKM expression, also because in the finite type 
case (i.e. in the $B_r$-representation rather than the $B_\infty$--one) many 
technical issues arise, at the point that in \cite{gln}, unlike in the DJKM 
case, the authors are not able to provide a formula for the generating 
function of the elementary endomorphisms tout-court, but only the generating 
function of all the images of a specified basis element through the 
elementary endomorphisms. This causes an unpleasant dependence of their 
formula from the partition parametrizing the basis element whose image is 
computed. Our idea is to remove the explicit dependence from the partition, 
by evaluating the image of the generating function of the Schur basis of $B_r
$. Again, we would have not able to achieve the goal  without heavily using 
the formalism of Schubert derivations as in \cite{BeGa,BeNa,gln,SDIWP}. 

We 
should finally remark that the way Schubert derivations remove the necessity 
to deal with partial derivatives,  carries a big potential to extend our main results to tropical situations  as indicated e.g. in 
\cite{GaRow1}. As a matter of fact, in \cite{GaRow1} a Grassmann semi-algebra is constructed 
within the framework of systems and the constructions lends itself to extend 
the Schubert derivations. In addition. in the last section of the recent 
preprint \cite{GaRoCh}, a semi-algebra version of our Theorem~
\ref{GattoSalehv1} is also provided,  essentially due to the fact that the techniques exposed in \cite{gln} work
in that more constrained situation. It seems likely that a sharpening of \cite[Theorem 7.24]{GaRoCh} can be naturally achieved using Cauchy type formulas for polynomial semi-algebras, a task which we temporarily postpone to further investigations.

\claim{\bf Structure of the paper.} We collect in Section~\ref{sec1:notation} 
the minimal background to follow the proofs of the main result, all based on 
the manipulation with Schubert derivations. Section~\ref{sec:sec3} is just a 
reformulation of the main formula in \cite{gln}. Section \ref{sec:sec4} 
contains the proof of the main theorem as well as the proof of many technical 
lemmas which on one hand are interesting in their own and, on the other hand, 
can be read within the classical theory of symmetric functions as in the 
classical reference \cite{MacDonald}. Section \ref{sec:sec5} is finally 
devoted to rephrase the DJKM representation with the purpose to present it 
within a unified perspective  together with Theorem~\ref{thm:mnthm1}.

\subsection*{Acknowledgments} We thank Letterio Gatto for having pointed out the papers \cite{GaRoCh, GaRow1}, so suggesting possible new research directions, and Parham Salehyan and Inna Scherbak for useful discussions and comments.

\section{Preliminaries and Notation}\label{sec1:notation}

\claim{} Let $X$ be an indeterminate over $\QQ$. We denote by $\Vcal:=\QQ[X^{-1},X]$ the vector space of Laurent polynomials, 
with basis $(X^i)_{i\in\ZZ}$ and for all $j\in\ZZ$ we write $\d^j$ for the 
unique linear form on $\Vcal$ such that $\d^j(X^i)=\delta^{ij}$. Let $V:=
\QQ[X]$ be the vector space of polynomials. It is a vector subspace of $\Vcal
$. The vector spaces $\Vcal^*:=\bigoplus_{j\in\ZZ}\QQ\, \d^j$ and $V^*:=
\bigoplus_{j\geq 0}\QQ \d^j$ are the {\em  restricted duals} of $\Vcal $ and 
$V$ respectively. Let $gl(V)$ and $gl(\Vcal)$ be respectively $V\otimes V^*$ and $\Vcal\otimes\Vcal^*$.

We denote by $\bfX(z)$ and $\bmd(w^{-1})$ the  {\em generating series} of the  basis 
elements of $V$ and of $V^*$ respectively, i.e.:
\begin{equation*}
\ubX(z):=\sum_{i\geq 0}X^iz^i\qquad \mathrm{and}\qquad \bmd(w^{-1}):=\sum_{j
\geq 0}
\d^j w^{-j}.\label{eq2:generating}
\end{equation*}
\claim{\bf Partitions and exterior algebras.} We denote by $\Pcal$  the set 
of all partitions. This is the set of all non-increasing sequences $
\blamb:=(\lambda_1\geq \lambda_2\geq\lambda_3\geq\cdots)$ of integers whose 
terms are all zero but finitely many. The terms of $\blamb$ are its {\em 
parts}. The length $\ell(\blamb)$ is the number of non zero parts and $
\Pcal_r$ stands for the set of all partitions with length at most $r$.
Let $\bwV=\bigoplus_{i\geq 0}\bw^iV$ be the {\em exterior algebra} of $V$. It 
is a graded algebra where $\bw^0V=\QQ$ and, for all $r\geq 1$:
$$
 \bw^rV:=\bigoplus_{\blamb\in
\Pcal_r}\QQ\cdot X^r(\blamb),
$$
where we have used the notation:
$$
X^r(\blamb):=X^{r-1+\lambda_1}\wedge X^{r-2+\lambda_2}\wedge\dots\wedge 
X^{\lambda_r}.
$$
In particular $\bw^1V=\bigoplus\QQ\cdot \bfX(\lambda)=\bigoplus\QQ\cdot X^
\lambda=V$.
\claim{\bf The ring $B_r$.}\label{TheringBr} For $r\geq 1$, let  $B_r:=
\QQ[e_1,
\ldots,e_r]$ be the polynomial ring in the $r$ indeterminates $(e_1,\ldots,e_r)
$. By convention one sets $B_0=\QQ$.
Consider the generic polynomial $E_r(z):=1-e_1z+\cdots+(-1)^r e_r z^r\in 
B_r[z]$, and the sequences $H_r:=(h_j)_{j\in\ZZ}$ and $\bfx_r=(x_i)_{i\in\ZZ}
$ defined through the equality:
\begin{equation*}
\sum_{n\in\ZZ}h_nz^n:={1\over E_r(z)}=\exp \left(\sum_{i \geq 1}x_i z^i
\right),\label{eq:defhj}
\end{equation*}
holding in $B_r\llb z\rrb.$
In particular $h_j=0$ if $j<0$ and $h_0=1$. Moreover for $j\geq 0$, $h_j$ is 
an 
explicit \emph{isobaric} polynomial of degree $j$ in $(e_1,\ldots, e_r)$, 
once one gives the weight  of $e_j$ to be $j$.

\claim{} It is well known that the {\em Schur determinants} 
\begin{equation*}
\Delta_\blamb(H_r):=\det(h_{\lambda_j-j+i})_{1\leq i,j\leq r}=\begin{vmatrix}
	h_{\lambda_1} & h_{\lambda_2-1} &\ldots & h_{\lambda_r -r+1} \\
	h_{\lambda_1+1} & h_{\lambda_2} &\ldots & h_{\lambda_r-r+2} \\
	\vdots  & \vdots & \ddots & \vdots \\
	
	h_{\lambda_1+r-1} & h_{\lambda_2+r-2} &\ldots &  h_{\lambda_r} \\
\end{vmatrix},\label{eq1:schurdet}
\end{equation*}
form a $\QQ$--basis of $B_r$ parametrized by 
the partitions of length at most $r$, i.e.:
\begin{equation*}
B_r:=\bigoplus_{\blamb\in\Pcal_r}\QQ\cdot \Delta_\blamb(H_r)\label{eq1:modBr}
\end{equation*}
 Thus the linear extension of the sets map 
\begin{equation*}
\Delta_\blamb(H_r)\mapsto X^r(\blamb).\label{eq0:fbfc}
\end{equation*}
gives a natural $\QQ$-vector space isomorphism $B_r\sra \bw^rV$.

If $r=\infty$, one sets
$$
B=B_\infty=\QQ[x_1,x_2,\ldots]
$$
and denotes by $S_\blamb(\bfx):=\mathrm{det}(S_{\lambda_j+j-1}(\bfx))$ the $\QQ$-basis element of $B$ corresponding to the partition $\blamb$, where
the sequence $(S_1(\bfx),S_2(\bfx),\dots)$ is obtained through the equation 
$$
\sum_{j\in\mathbb{Z}}S_j(\bfx)z^j=\mathrm{exp}\left(\sum_{i\geq 1}x_iz^i\right).
$$
\claim{\bf Schur Polynomials.}\label{ExtSchurPolys} Let $\bfz_k:=(z_1,\ldots, z_k)$ be an ordered finite sequence of formal variables and consider
the $\bw^k V$-valued formal power series
$$
\ubX(z_k)\w\cdots\w \ubX(z_1).
$$
It vanishes along all the diagonals $z_i-z_j=0$ ($i\neq j$). Therefore it is 
divisible by the Vandermonde determinant $\Delta_0(\bfz_k)=\prod_{1\leq 
	i<j\leq k}(z_j-z_i)$. 
The equality 

\begin{equation*}
\sum_{\bmu\in \Pcal_k} X^k(\bmu)s_\blamb(\bfz_k)\Delta_0(\bfz_k)=\ubX(z_k)\w\cdots\w \ubX(z_1),
\label{eq1:gfbz}
\end{equation*}
define the {\em Schur symmetric polynomial} $s_\blamb(\bfx_k)$. This definition coincides with the usual one as in \cite[formula (3.1)]{MacDonald}.

\smallskip
%
%

\claim{\bf The $\End_\QQ(V)$-module structure of $\bwV$.} \label{sec:emstru} 
Given $\phi\in \End_\QQ(V)$ let us denote by $\delta(\phi)$ the unique derivation of $\bwV$ such that $\delta(\phi)_{|_V}=\phi$ (see \cite[Section 3.1]{BeCoGaVi}). In other words, $\delta(\phi)(u\w v)=\delta(\phi)u\w v+u\w\delta(\phi)v$ for all $u,v\in\bwV$ and $\delta(\phi)w=\phi(w)$ for all $w\in V=\bw^1V$.
\bclm{\bf Proposition} {\em The plethystic exponential of $\delta(\phi)\in \End_\QQ(\bwV)$
$$
\Dcal^\phi(z):={\mathrm Exp}(\delta(\phi)z):=\exp\left(\sum_{i\geq 1}{1\over i}\delta(\phi^i)z^i\right):\bwV\sra \bwV\llb z\rrb
$$
is the unique Hasse-Schmidt derivation on $\bwV$ such that 
$$
\Dcal^\phi(z)_{|V}=\sum_{i\geq 0}\phi^iz^i\in \End_\QQ(V).
$$
}
\eclm
See \cite{BeThesis} for the proof. Recall by e.g. \cite{SCHSD} or \cite{HSDGA}  that to say $\Dcal^\phi(z)$ is a Hasse-Schmidt (HS) derivation means that
$$
\Dcal^\phi(z)(u\w v)=\Dcal^\phi(z)u\w \Dcal^\phi(z)v\qquad \forall u,v\in\bwV
$$
In particular, putting $\overline{\Dcal}^\phi(z)={\exp}(-\delta(\phi)z)$, which is easily seen to be a HS derivation as well, the integration by parts formula holds:
\be
\overline{\Dcal}^\phi(z)u\w v=\Dcal^\phi(z)(u\w \overline{\Dcal}^\phi(z)v)
\label{eq:intp1}
\ee

\claim{} By abuse of notation we denote by $X$ and $X^{-1}$ the $\QQ$-linear maps $V\sra V$ given by multiplication by $X$ and $X^{-1}$ respectively, where for all $i>0$
\smallskip
$$
X^{-i}X^j=\left\{\begin{matrix}X^{j-i}& if& i\leq j\cr\cr 0&if &i>j
\end{matrix}\right.
$$
Consider
\begin{eqnarray}
		\sigma_+(z)&=&\sum_{i\geq 0}\sigma_iz^i:=\exp\left(\sum_{i\geq 1}{1\over i}
		\delta(X^i)z^i\right):\bwV\sra\bwV\llb z\rrb,  \\  \sigma_-(z)&=&\sum_{i\geq 0}\sigma_{-i}z^{-i}:=
		\exp\left(\sum_{i\geq 1}{1\over i}\delta(X^{-i})z^{-i}\right):\bwV\sra\bwV[ z^{-1}],
\end{eqnarray}
		and their inverses as elements of $\End_\QQ(\bwV)\llb z^{\pm 1}\rrb$
\begin{eqnarray}
		\ovsig_+(z)&=&\sum_{i\geq 0}(-1)^i\ovsig_iz^i:=\exp\left(-\sum_{i\geq 1}{1\over i}
		\delta(X^i)z^i\right):\bwV\sra\bwV [z],  \\  \ovsig_-(z)&=&\sum_{i\geq 0}(-1)^i\ovsig_{-i}z^{-i}:=
		\exp\left(-\sum_{i\geq 1}{1\over i}\delta(X^{-i})z^{-i}\right):\bwV\sra\bwV[ z^{-1}].
\end{eqnarray}
		
\begin{claim}{\bf Proposition}(\cite[Proposition 2.2]{pluckercone}). {\em The maps $\sigma_{\pm}(z)$ and $\ovsig_{\pm}(z)$ are the unique  (HS) derivations on the exterior algebra $\bwV$ such that 
	\be
	\sigma_+(z)X^j=\sum_{i\geq 0}X^{j+i}z^i, \qquad\qquad \ovsig_+(z)X^j=X^j-X^{j+1}z,\label{eq0:s+bi}
	\ee
	and
	\be
	\sigma_-(z)X^j=\sum_{i\geq 0} X^{j-i} z^{-i}, \qquad\qquad  \ovsig_-(z)X^j=X^j- X^{j-1}z^{-1},\label{eq0:s-bi}
	\ee
	putting $X^i=0$ for $i<0$. 
They are called {\em Schubert derivations}.
}
\eclm

\smallskip

%
%
	\claim{\bf Transposition.} \label{sec:transp}The transpose  $\sigma_\pm(z)^T:\bwV^*\sra\bwV^*\llb z^{\pm 1}\rrb$ of the Schubert derivation $
	\sigma_\pm(z)$ is defined via its action on homogeneous elements. If $\eta\in \bw^rV^*$, 
	then one stipulates that $\sigma_\pm(z)^T\eta(u)=\eta(\sigma_\pm(z)u)$, for all $u\in \bw^r\hskip-2pt V$. By 
	\cite[Proposition 2.8]{pluckercone} $\sigma_\pm(z)^T$ is a HS--derivation of $\bwV^*$.	In the sequel we will need the fact that
	$$\sigma_-(z)^TX^j=\sum_{i\geq 0}X^{j+i}z^{-i}.$$
	\claim{\bf The  $B_r$-module structure of $\bw^rV$.}\label{sec3:mosst} For all $\bfu\in\bw^rV$, define 
	\be
	e_i\bfu=\ovsig_i\bfu\quad  \mathrm{or, \,\,equivalently,}\quad h_i\bfu=\sigma_i\bfu.\label{eq:prems}
	\ee
	In particular:
	$$
	\ovsig_+(z)\bfu={E_r(z)}\cdot \bfu\qquad\mathrm{and}\qquad \sigma_+(z)\bfu:={1\over E_r(z)}\bfu,
	\qquad \forall \bfu\in\bw^rV .
	$$

	\bclm{\bf Proposition} \label{prop:giamb}. {\em Equations \eqref{eq:prems} make $\bw^rV$ into a free $B_r$-module generated by $X^r(0):=X^{r-1}\w\cdots\w X^0$.
	}
	\eclm
\proof It is a consequence of the fact, explained in \cite{SCHSD}, see also \cite[Proposition 3.5]{BeCoGaVi}, that
	Giambelli's formula for the Schubert 
		derivation $\sigma_+(z)$ holds:
		\be
		X^r(\blamb)=\Delta_\blamb(\sigma_+(z))\cdot X^r(0):=\left(\det(\sigma_{\lambda_j-j+i})_{1\leq i,j\leq r}\right)\cdot X^r(0).
		\label{eq:gmbsh}
		\ee
		Therefore 
		$$
		h_i\cdot X^r(\blamb)=\sigma_iX^r(\blamb)=\sum_{\bmu\in\mathcal{P}_r} X^r(\bmu)=\sum_{\bmu\in\mathcal{P}_r} \Delta_\bmu(H_r)X^r(0)=(h_i\Delta_\blamb(H_r))X^r(0)
		$$
		Hence $\bw^rV$ is a free $B_r$-module of rank $1$ generated by $	X^r(0)$.\qed

	The fact that $\bw^rV$ is a free $B_r$-module of rank $1$ generated by $X^r(0)$, as prescribed by equality \eqref{eq:gmbsh}, shows that the Schubert derivations  $\sigma_-(z), \ovsig_-(z)$ induce maps $B_r\sra B_r[z^{-1}]$ that, abusing notation, will be denoted in the same way. Their action on a basis element $\Delta_\blamb(H_r)$ of $B_r$ is defined through its action on $\bw^rV$:
	\begin{eqnarray}
	(\ovsig_-(z)\Delta_\blamb(H_r))X^r(0)&:=&\ovsig_-(z)X^r(\blamb),\label{eq3:ovsig-b}\\ 
	\cr (\sigma_-(z)\Delta_\blamb(H_r))X^r(0)&:=&\sigma_-(z)X^r(\blamb).\label{eq3:sig-b}
	\end{eqnarray}
	Denote by $\ovsig_-(z)H_r$ (respectively $\sigma_-(z)H_r$)  the sequence  $(\ovsig_-(z)h_j)_{j\in\ZZ}$ (respectively\linebreak  $(\sigma_-(z)h_j)_{j\in\ZZ}$). By using \cite[Theorem 5.7]{pluckercone}, and exploiting the Laksov \& Thorup determinantal formula as in \cite[Main Theorem 0.1]{LakTh01}, one obtains the  following statement, which gives a practical way to evaluate the image of $\Delta_\blamb(H_r)$ through the maps $\ovsig_-(z)$ and $\sigma_-(z)$ defined by \eqref{eq3:sig-b} and \eqref{eq3:ovsig-b}.
	
	\claim{\bf Proposition} (\cite[Proposition 5.3]{pluckercone}){\bf .} \label{prop3:prop35} \label{prop3:prop35} {\em 
		For all $r\geq 0$ and all $\blamb\in\Pcal_r$
		$$
		\sigma_-(z)h_j=\sum_{i\geq 0}h_{j-i}z^{-i}\qquad \mathrm{and}\qquad \ovsig_-(z)h_j=h_j-{h_{j-1}z^{-1}}.
		$$
		Moreover,
		$$
		\sigma_-(z)\Delta_\blamb(H_r)=\Delta_\blamb(\sigma_-(z)H_r)\qquad\mathrm{and}\qquad \ovsig_-(z)\Delta_\blamb(H_r)=\Delta_\blamb(\ovsig_-(z)H_r).\label{eq1:comm1}
		$$
	}

\section{The $gl(V)$-structure of $B_r$ revisited 1}\label{sec:sec3}

In this section we revisit \cite[Theorems 5.7 and 6.4]{gln} and \cite[Section 9]{SDIWP}, getting a more explicitly and elegant expression
to realize the $gl(V):=V\otimes V^*$-module structure of the ring $B_r$.

 Let $\beta\in V^*$. In the paper \cite[Section 5.1 and Lemma 5.3]{BeCoGaVi} one learns to phrase the usual contraction endomorphism
 $$\beta\lrcorner:\bigwedge V\rightarrow\bigwedge V$$ 
 via  the following diagram
\be\label{contraction1}
\beta\lrcorner\ubX^{r}(\blamb):=\left\vert
\begin{array}{cccc}
\beta(X^{r-1+\lambda_1}) & \beta(X^{r-2+\lambda_2})  & \cdots & \beta(X^{\lambda_r}) \cr\cr
X^{r-1+\lambda_1} & X^{r-2+\lambda_2}  & \cdots & X^{\lambda_r} 
\end{array}
\right\vert ,
\ee
which means that the scalar $(-1)^{i+1}\beta(X^{r-i+\lambda_i})$ is the coefficient of the element of $\bigwedge^{r-1}V$ obtained by removing
the $i$-th  exterior factor of the wedge product of the elements
in the second row, namely $X^{r-1+\lambda_1}\wedge\dots\wedge X^{\lambda_r}=X^{r}(\blamb)$. 


For example, it follows by the very definition \eqref{contraction1} that $$\d^0\lrcorner X^r(0)=(-1)^{r-1} X^{r-1}((1^{r-1})).$$


\claim{} The Lie algebra $gl(V)$ acts on the exterior algebra via the map $\delta$ as in Section \ref{sec:emstru}. In particular
$$
\delta(X^i\otimes \d^j)(u)=X^i\w \d^j\lrcorner u\qquad \forall u\in\bwV .
$$
Using the definition of the Schubert derivation $\sigma_+(z)$ as in \eqref{eq0:s+bi}, the generating function $\Ecal(z,w^{-1})$ of the basis $(X^i\otimes \d^j)_{i,j\geq 0}$ of $gl(V)$ can be written as
$$
\mathcal{E}(z,w^{-1})=\bfX(z)\otimes\bmd(w^{-1})=\sigma_+(z)X^0\otimes\bmd(w^{-1})
$$
and acts on $\bw V$ as
\be\label{glmodulebr}
\Ecal(z,w^{-1})(X^r(\blamb))=\sigma_+(z)X^0\w (\bmd(w^{-1})\lrcorner X^r(\blamb)).
\ee

\medskip

\noindent The following rephrases in a more elegant and transparent way the description of the $gl(V)$ structure of $\bw^rV$ proposed in \cite[Theorem ~4.3]{{gln}}.

\claim{\bf Proposition.}\label{GattoSalehv1}
{\em The action of $\Ecal(z,w^{-1})$ on the basis element $X^r(\blamb)$ is given by:
\begin{equation*}
    \mathcal{E}(z,w^{-1})X^r(\blamb)=\frac{z^{r-1}}{w^{r-1}}\sigma_+(z)\overline{\sigma}_-(z)\left(\begin{vmatrix}
            w^{-\lambda_1} & w^{-\lambda_2+1} &\dots & w^{r-1-\lambda_r} & 0\cr\cr
            X^{r+\lambda_1} & X^{r-1+\lambda_2} & \dots & X^{1+\lambda_r} & X^0
        \end{vmatrix}\right).
\end{equation*}

\begin{proof}
By applying equation \eqref{glmodulebr} and diagram \eqref{contraction1}, we obtain
\begin{eqnarray*}
        \mathcal{E}(z,w^{-1})X^r(\blamb) &=& \sigma_+(z)X^0\wedge (\bm{\d}(w^{-1})\lrcorner X^r(\blamb)) \cr\cr
        &=&\sigma_+(z)X^0\wedge \begin{vmatrix}
            \bmd(w^{-1})X^{r-1+\lambda_1} & \dots &  \bmd(w^{-1})X^{\lambda_r}\cr\cr
            X^{r-1+\lambda_1} & \dots & X^{\lambda_r}
        \end{vmatrix}
        \cr\cr\cr
        &=&\sigma_+(z)X^0\wedge \begin{vmatrix}
            w^{-r+1-\lambda_1} & \dots & w^{-\lambda_r}\cr\cr
            X^{r-1+\lambda_1} & \dots & X^{\lambda_r}
        \end{vmatrix}.
\end{eqnarray*}

\noindent Using the integration by parts \eqref{eq:intp1}  for the Schubert derivation $\sigma_+(z)$, one obtains:
\begin{eqnarray}
    \mathcal{E}(z,w^{-1})X^r(\blamb) 
        &=&\sigma_+(z)\left(X^0\wedge \overline{\sigma}_+(z)\begin{vmatrix}
            w^{-r+1-\lambda_1} & \dots & w^{-\lambda_r}\cr\cr
            X^{r-1+\lambda_1} & \dots & X^{\lambda_r}
        \end{vmatrix}\right)\cr\cr
        &=&\sigma_+(z)\left((-1)^{r-1}\begin{vmatrix}
            w^{-r+1-\lambda_1} & \dots & w^{-\lambda_r}\cr\cr
            \overline{\sigma}_+(z)X^{r-1+\lambda_1} & \dots & \overline{\sigma}_+(z)X^{\lambda_r}
        \end{vmatrix}\wedge X^0\right).\label{eq:lastd}
            \end{eqnarray}
By \eqref{eq0:s+bi} and \eqref{eq0:s-bi} one has $\overline\sigma_{+}(z)X^{j}=-z\overline\sigma_{-}(z)X^{j+1}$. Hence \eqref{eq:lastd}
 can be rewritten as
\begin{equation*}
    \begin{split}
        &\sigma_+(z)\left((-1)^{r-1} \begin{vmatrix}
            w^{-r+1-\lambda_1} & \dots & w^{-\lambda_r}\cr\cr
            -z\overline{\sigma}_-(z)X^{r+\lambda_1} & \dots & -z\overline{\sigma}_-(z)X^{1+\lambda_r}
        \end{vmatrix}\wedge X^0\right)\cr\cr
        &=\sigma_+(z)\left(z^{r-1} \begin{vmatrix}
            w^{-r+1-\lambda_1} & \dots & w^{-\lambda_r} & 0\cr\cr
            \overline{\sigma}_-(z)X^{r+\lambda_1} & \dots & \overline{\sigma}_-(z)X^{1+\lambda_r} & X^0
        \end{vmatrix}\right)\cr\cr
        &=\frac{z^{r-1}}{w^{r-1}}\sigma_+(z)\overline{\sigma}_-(z)\begin{vmatrix}
            w^{-\lambda_1} & w^{-\lambda_2+1} &\dots & w^{r-1-\lambda_r} & 0\cr\cr
            X^{r+\lambda_1} & X^{r-1+\lambda_2} & \dots & X^{1+\lambda_r} & X^0
        \end{vmatrix}
    \end{split}
\end{equation*}
which is the desired expression.\qed
\end{proof}
\claim{\bf Example.} {\em Applying Proposition \ref{GattoSalehv1} to the particular case $r=2$, and using 
the properties of Schubert derivations in equations \eqref{eq0:s+bi} and \eqref{eq0:s-bi}, one gets 
\begin{eqnarray*}
\Ecal(z,w^{-1})X^2(\blamb)& = &\frac{z}{w}\sigma_+(z)\overline{\sigma}_-(z)\left(\begin{vmatrix}
            w^{-\lambda_1} & w^{-\lambda_2+1} & 0\cr\cr
            X^{2+\lambda_1} & X^{1+\lambda_2} & X^0
        \end{vmatrix}\right) \\
        &=&\frac{z}{w}\sigma_+(z)\overline{\sigma}_-(z)\big(w^{-\lambda_1}(X^{1+\lambda_2}\wedge X^0)-w^{-\lambda_2+1}(X^{2+\lambda_1}\wedge X^0)\big) \\
        &=& \frac{z}{w^{1+\lambda_1}}\big(\sigma_+(z)(X^{1+\lambda_2}-X^{\lambda_2}z^{-1})\wedge \sigma_+(z)X^0\big) \\
        & &- \frac{z}{w^{\lambda_2}}\big(\sigma_+(z)(X^{2+\lambda_1}-X^{1+\lambda_1}z^{-1})\wedge\sigma_+(z)X^0\big)  \\
        &=& \frac{-1}{w^{1+\lambda_1}}\big(X^{\lambda_2}\wedge \sigma_+(z)X^0\big)+\frac{1}{w^{\lambda_2}}\big(X^{1+\lambda_1}\wedge\sigma_+(z)X^0\big)  \\
        &=& \frac{-1}{w^{1+\lambda_1}}\big(X^{\lambda_2}\wedge \sigma_+(z)X^0\big)+\frac{1}{w^{\lambda_2}}\big(X^{1+\lambda_1}\wedge\sigma_+(z)X^0\big)  \\
         &=& \left(\frac{X^{1+\lambda_1}}{w^{\lambda_2}}-\frac{X^{\lambda_2}}{w^{1+\lambda_1}}\right)\wedge \sigma_+(z)X^0 \\
         &=& \left(\frac{X^{1+\lambda_1}}{w^{\lambda_2}}-\frac{X^{\lambda_2}}{w^{1+\lambda_1}}\right)\wedge \bfX(z).
\end{eqnarray*}
}

\claim{\bf Remark.}
Besides the formula displayed in Proposition \ref{GattoSalehv1} being very manageable and explicit, 
there is still a dependence on $\bm{\lambda}$ that is desirable to remove. The idea to do this is to sum all the above expression for $\blamb$ ranging over all the partitions. Better, we will apply the generating series $\Ecal(z,w)$ to a the generating function $\sum_{\blamb\in\Pcal}X^r(\blamb)s_\blamb(\bft_r)$. By \cite[Lemma ~8.2]{{BeCoGaVi}} (which is essentially an exponential way to write the Cauchy formula as in \cite[Formula (3), p.~52]{Fulyoung}) that is
$$
\exp\left(\sum_{i\geq 1}x_ip_i(\bft_r)\right)X^r(0),
$$
where $p_i(\bft_r):=t_1^i+\cdots+t_r^i$.

\section{The  $gl(V)$ structure of $B_r$ revisited 2}\label{sec:sec4}
Recall that the $B_r$-module structure  of $\bw^rV$, established in Proposition~\ref{prop:giamb}, maps $B_r$--valued formal power series to $\bw^rV$-valued formal power series via the map $\phi\mapsto \phi\cdot X^r(0)$.

\bclm{\bf Definition.}\label{def:def41} {\em We denote by $\Ecal_r(z,w^{-1},\bft_r)$ the unique $B_r$-valued formal power series in the indeterminates $(z,w^{-1},\bft_r)$ such that
\begin{eqnarray*}
\Ecal_r(z,w^{-1},\bft_r)X^r(0)&=&\sum_{((i,j),\blamb)\in\NN^2\times\Pcal_r}\left[X^i\otimes \d^j\star \Delta_\blamb(H_r)\right]z^iw^{-j}s_\blamb(\bft_r)\bfX^{r}(0) \cr\cr\cr 
&=&\sum_{\blamb\in\Pcal_r} \bfX(z)\w \big(\bmd(w^{-1})\lrcorner \bfX^r(\blamb)s_\blamb(\bft_r)\big).
\end{eqnarray*}
}
\eclm

\noindent Clearly $\Ecal(z,w^{-1},\bft_r)\in B_r\llb z, w^{-1},\bft_r\rrb$ and the purpose of this section  is to determine an exponential-like expression for
$$\Ecal_r(z,w^{-1},\bft_r):=\sum_{((i,j),\blamb)\in\NN^2\times\Pcal_r}\left[X^i\otimes \d^j\star \Delta_\blamb(H_r)\right]z^iw^{-j}s_\blamb(\bft_r)\in B_r\llb z,w^{-1},\bft_r\rrb.
$$

%
%

\noindent To this purpose we first invoke \cite[Lemma ~8.2]{{BeCoGaVi}} according which
\be\label{invoke}
\sum_{\blamb\in\Pcal_r}\bfX^r(\blamb)s_\blamb(\bft_r)=\sigma_+(\bft_r)\bfX^r(0).
\ee
Recall that the symbols $\sigma_{\pm}(\bft_r)$ and  $\overline{\sigma}_{\pm}(\bft_r)$ stand for \emph{multivariables Schubert derivations}  (see \cite[Section 4.1]{BeCoGaVi}):
\begin{eqnarray*}
\sigma_{\pm}(\bft_r):=\sigma_{\pm}(t_1)\cdots\sigma_{\pm}(t_r) & \ \ \mbox{and} \ \ &\overline{\sigma}_{\pm}(\bft_r):=\overline{\sigma}_{\pm}(t_1)\cdots\overline{\sigma}_{\pm}(t_r).
\end{eqnarray*}
They are HS derivations on $\bigwedge V$, in the sense that $\sigma_{\pm}(\bft_r)(u\wedge v)=\sigma_{\pm}(\bft_r)(u)\wedge \sigma_{\pm}(\bft_r)(v)$, holding the same to their inverses
$\overline{\sigma}{\pm}(\bft_r)$.

Next we need some preliminaries results.


\claim{\bf Lemma. }{\em 
Let $t$ be one formal variable over $\mathbb{Q}$. The following equality holds for all $r\geq 2$:
\begin{equation*}
    \sigma_+(t)\bfX^r(0)=\sigma_+(t)X^{r-1}\wedge {\bfX}^{r-1}(0).
\end{equation*}
}

\begin{proof}
One argues by induction on $r$. The case  $r=1$ is ruled by the equality
\begin{equation}\label{eq:sigma+X^i}
    \sigma_+(t)X^i=X^i+t\sigma_+(t)X^{i+1}.
\end{equation}
Now suppose that the property holds for $r-1\geq 1$, then:
\begin{center}
\begin{tabular}{rllrr}
$\sigma_+(t)\bfX^r(0)$&$=$&$\sigma_+(t)X^{r-1}\wedge \sigma_+(t) \bfX^{r-1}(0)$ \ \ \  (because $\sigma_+(t)\in HS(\bwV))$\cr\cr
&$=$&$\sigma_+(t)X^{r-1}\wedge(\sigma_+(t)X^{r-2}\wedge \bfX^{r-2}(0))$ \ \ \  (induction hypothesis) &\cr\cr
&$=$&$\sigma_+(t)X^{r-1}\wedge(\sigma_+(t)X^{r-2}\wedge X^{r-3}\wedge \dots\wedge X^0)$\cr\cr
&$=$&$\sigma_+(t)X^{r-1}\wedge(X^{r-2}+t\sigma_+(t)X^{r-1})\wedge X^{r-3}\wedge \dots\wedge X^0$\ \ \ (eq. (\ref{eq:sigma+X^i}))& \cr\cr
&$=$&$\sigma_+(t)X^{r-1}\wedge \bfX^{r-1}(0)$
\end{tabular}
\end{center}
as desired.\qed
\end{proof}


  
\claim{\bf Proposition. }\label{name_4}
{\em The multivariable Schubert derivations satisfy the following equalities:
\begin{eqnarray}
 \ovsig_+(\bft_r)X^0&=&X^0+\displaystyle{\sum_{i=1}}^r(-1)^ie_i(\bft_r)X^i  \label{eq1:sim1}\\ \cr
 \sigma_+(\bft_r)X^0&=&X^0+\displaystyle{\sum_{i\geq 1}}h_i(\bft_r)X^i, \label{eq2:sim2}
\end{eqnarray}
where $e_i(\bft_r)$ and $h_i(\bft_r)$ are, respectively, the elementary and complete symmetric polynomial of degree $i$ in the indeterminates $t_1,\ldots,t_r$.
}
\medskip

\begin{proof}  For $r=1$ is the definition of $\sigma_+(t)$ and $\ovsig_+(t)$. The general case follows by induction on $r$.

Indeed, if \eqref{eq1:sim1} holds for $r-1\geq 1$,  then
\begin{equation*}
    \begin{split}
        \overline{\sigma}_+(\mathbf{t}_r)X^0&=\overline{\sigma}_+(t_r)\overline{\sigma}_+(\mathbf{t}_{r-1})X^0= \sum_{i=0}^{r-1}(-1)^i(X^i-t_rX^{i+1})e_i(\mathbf{t}_{r-1})\\
        &=\sum_{i=0}^{r-1}(-1)^ie_i(\mathbf{t}_{r-1})X^i+\sum_{i=0}^{r-1} (-1)^{i+1}t_re_i(\mathbf{t}_{r-1})X^{i+1}\\
        &=\sum_{i=0}^{r-1}(-1)^ie_i(\mathbf{t}_{r-1})X^i+\sum_{i=1}^{r} (-1)^it_re_{i-1}(\mathbf{t}_{r-1})X^i\\
        &=X^0+\sum_{i=1}^{r-1}(-1)^i \left(e_i(\mathbf{t}_{r-1})+t_re_{i-1}(\mathbf{t}_{r-1})\right)X^i+(-1)^rt_re_{r-1}(\mathbf{t}_{r-1})X^r\\
        &=\sum_{i=0}^{r} (-1)^iX^ie_i(\mathbf{t}_r).
    \end{split}
\end{equation*}
Equality \eqref{eq2:sim2} is a consequence of \eqref{eq1:sim1}, because $\sigma_+(\bft_r)$ and $\ovsig_+(\bft_r)$ are mutually inverses  in $\End_\QQ(\bwV)\llb \bft_r\rrb$.\qed
\end{proof}

\bclm{\bf Lemma.}\label{lem:lemel} {\em In the ring of formal power series $\QQ\llb u\rrb$ the following equality holds:
$$
1-{u}=\exp\left(-\sum_{i\geq 1}{u^i\over i}\right).\label{eq:explog}
$$
}
\eclm
\proof 
This is well known. Just write
\begin{eqnarray*}
1-u=\exp\left(\log(1-u)\right)=\exp\left(-\int{du\over 1-u}\right)
&=&\exp\left(-(1+u+u^2+\cdots)du\right)\cr\cr
&=&\exp\left(-\sum_{i\geq 1}{u^i\over i}\right).\hskip 100pt\qed
\end{eqnarray*}
\medskip

\claim{\bf Proposition. }\label{name_1}
The following equalities hold:

\be \sigma_-(w)\sigma_+(\mathbf{t}_r)\bfX^r(0)=\exp\left(\sum_{n\geq 1}\frac{1}{n}p_n\left(\frac{\mathbf{t}_r}{w}\right)\right)\sigma_+(\mathbf{t}_r)\bfX^r(0),\label{eq:sig-sig+(t_r)}
\ee
\be\overline{\sigma}_-(w)\sigma_+(\mathbf{t}_r)\bfX^r(0)=\exp\left(-\sum_{n\geq 1}\frac{1}{n}p_n\left(\frac{\mathbf{t}_r}{w}\right)\right)\sigma_+(\mathbf{t}_r)\bfX^r(0).\label{eq:sig-sig+(t_r)2}
\ee
\medskip
\proof
We first prove that for a single formal variable $t$ one has:
$$
\sigma_-(w)\sigma_+(t)X^0=\exp\left(\sum_{i\geq 1}{t^i\over w^i}\right)\sigma_+(t)X^0
$$
Indeed by the very definition of $\sigma_+(t)$:
\begin{eqnarray*}
\sigma_-(w)\sigma_+(t)X^0&=&\sigma_-(w)\left(\sum_{i\geq 0}X^{i}t^i\right) =X^0+\sigma_-(w)\left(\sum_{i\geq 0}X^{i+1}t^{i+1}\right)\cr
&=&\sum_{i\geq 0}\sum_{j=0}^{i}X^{i-j}{t^i\over w^j}=\left(\sum_{i\geq 0}{t^i\over w^i}\right)\sigma_+(t)X^0\cr\cr
&=&{1\over 1-\displaystyle{t\over w}}\sigma_+(t)X^0=\exp\left(\sum_{i\geq 1}{t^i\over w^i}\right)\sigma_+(t)X^0\label{eq:r=1}
\end{eqnarray*}
where in the last equality we used Lemma~\ref{lem:lemel}.
For $r>1$, denote by
 $\Delta_0(\bft_r)$ the Vandermonde determinant
\be
\Delta_0(\bft_r)=\prod_{1\leq i<j\leq r}(t_j-t_i),\label{eq:vande}
\ee

 we have:
\begin{eqnarray*}
    \Delta_0(\bft_r)\sigma_-(w)\sigma_+(\bft_r)X^r(0)&=&\sigma_-(w)(\sigma_+(t_1)X^0\w\dots\w\sigma_+(t_r)X^0)\cr\cr&=&\sigma_-(w)\sigma_+(t_1)X^{0} \w \cdots \w  \sigma_-(w)\sigma_+(t_r)X^{0}
    \end{eqnarray*}

and because  $\sigma_-(w)$ is a HS-derivation we have
\begin{eqnarray*}
         \Delta_0(\bft_r)\sigma_-(w)\sigma_+(\mathbf{t}_r)X^r(0)&=&\sigma_-(w)\left(\sigma_+(t_1)X^{0}\w \cdots \w \sigma_+(t_r) X^0 \right)\cr\cr
         &=&\sigma_-(w)\sigma_+(t_1)X^{0} \w \cdots \w  \sigma_-(w)\sigma_+(t_r)X^{0}.
  \end{eqnarray*}
Now by virtue of Lemma \ref{lem:lemel}, the last above expression can be rewritten as 
\begin{equation*}
    \begin{split}
       &=\exp\left(\sum_{n\geq 1}{1 \over n}{t_1 ^n \over w^n}\right) \cdots \exp\left(\sum_{n\geq 1}{1 \over n}{t_r ^n \over w^n}\right) 
        \left(\sigma_+(t_1)X^{0}\w \cdots \w \sigma_+(t_r) X^0 \right)\\
        &=\Delta_0(\bft_r)\exp\left(\sum_{n\geq 1}{1 \over n}p_n\left({\bft_r \over w}\right)\right)\sigma_+(\bft_r)X^r(0).
    \end{split}
\end{equation*}
whence \eqref{eq:sig-sig+(t_r)} because $\Delta_0(\bft_r)$ is not a zero divisor in $\QQ[\bft_r]$.  The proof of \eqref{eq:sig-sig+(t_r)2} is similar and we omit it.  \qed


\claim{\bf Lemma. }\label{lemconstpart} {\em Let $\ell$ be a non-negative integer and set $(\ell^r)$ to be the partition with $r$ parts equal to $\ell$.
For each $r\geq 1$ we have the following identities:
\
\begin{enumerate}
\item\label{const_partition}$\overline{\sigma}_+(z)X^r(\ell^r)=X^r(\ell^r)+\sum_{i=1}^r (-1)^iz^i\,X^r((\ell+1)^i\ell^{r-i})$ and
\item $e_i X^r(\ell^r):=\overline{\sigma}_iX^r(\ell^r)=X^r((\ell+1)^i\ell^{r-i})$, for $0\leq i\leq r.$
\end{enumerate}
}

\medskip

\proof
The second item is a consequence of the first one. Assuming $r=1$,  item \ref{const_partition}) can be  seen as follows
\begin{equation*}
        \overline{\sigma}_+(z)X^1(\ell)=\overline{\sigma}_+(z)X^\ell=X^\ell-zX^{\ell+1}=X^1(\ell)-zX^1(\ell+1).
\end{equation*}
Assuming that  the statement holds true for some $r\geq 1$, the proof follows by induction on $r$. 

\medskip
\noindent $\overline{\sigma}_+(z)X^{r+1}(\ell^{r+1})$
\begin{equation*}
    \begin{split}
        &=\overline{\sigma}_+(z)(X^{\ell+r}\wedge X^r(\ell^r))\\ 
        &=(X^{\ell+r}-zX^{\ell+r+1})\wedge\overline{\sigma}_+(z)X^r(\ell^r) \hspace{3.5cm} \ (\overline{\sigma}_+(z)\in HS(\wedge V)\mbox { and eq.  } \eqref{eq0:s+bi}) \\ 
        &=(X^{\ell+r}-zX^{\ell+r+1})\wedge\left(X^r(\ell^r)+\sum_{i=1}^r (-1)^iz^i\,X^r(({\ell+1})^i\ell^{r-i})\right) \ \ \mbox{(induction hypothesis)}\\
        &=X^{\ell+r}\wedge X^r(\ell^r)-zX^{\ell+r+1}\wedge X^r(\ell^r)+\sum_{i=1}^r (-1)^{i+1}z^{i+1}\cdot X^{\ell+r+1}\wedge X^r(({\ell+1})^i\ell^{r-i})\\
        &=X^{r+1}(\ell^{r+1})-zX^{r+1}((\ell+1)\ell^r)+\sum_{i=2}^{r+1} (-1)^iz^i\cdot X^{\ell+r+1}\wedge X^r(({\ell+1})^{i-1}\ell^{r-(i-1)}) \\
        &=X^{r+1}(\ell^{r+1})-zX^{r+1}((\ell+1)\ell^r)+\sum_{i=2}^{r+1} (-1)^iz^i\cdot X^{r+1}(({\ell+1})^i\ell^{r+1-i})\\
        &=X^{r+1}(\ell^{r+1})+\sum_{i=1}^{r+1} (-1)^iz^i\cdot X^{r+1}(({\ell+1})^i\ell^{r+1-i})\hskip190pt \qed
    \end{split}
\end{equation*}


\claim{\bf Lemma. }\label{name_2}{\em The following commutations rule hold 
\begin{enumerate}
\item\label{d0_against_sigma} $\d^0\lrcorner \sigma_+(t)X^r(\blamb)=\sigma_+(t)(\d^0\lrcorner X^r(\blamb))$;
\item\label{d0_against_sigma_bar}
    $\d^0\lrcorner \overline{\sigma}_+(t)X^r(\blamb)=\overline{\sigma}_+(t)(\d^0\lrcorner X^r(\blamb))$.
\end{enumerate}
}

\begin{proof}
If the length of the partition is exactly $r$,  $l(\bm{\lambda})=r$, then both sides of (\ref{d0_against_sigma}) are equal to zero. If $l(\bm{\lambda})<r$, since  
$X^r(\blamb)=X^{r-1} (\blamb +(1^r))\w X^0$, then the left hand side of (\ref{d0_against_sigma}) is
\
\begin{equation*}
    \begin{split}
        \d^0\lrcorner\sigma_+(t)(X^{r-1}(\bm{\lambda}+(1^r))\wedge X^0)&=\d^0\lrcorner(\sigma_+(t)X^{r-1}(\bm{\lambda}+(1^r))\wedge \sigma_+(t)X^0)\cr\cr
%
&=(-1)^{r-1}\sigma_+(t)X^{r-1}(\bm{\lambda}+(1^r))\cr\cr
        &=\sigma_+(t)(\d^0\lrcorner X^r(\bm{\lambda})).\cr
    \end{split}
\end{equation*}
The proof of item \eqref{d0_against_sigma_bar} is completely analogous. \qed
\end{proof}

\claim{\bf Lemma. }\label{name} {\em It follows from \cite[Lemma ~5.5]{{BeCoGaVi}} that
\begin{equation*}
    \bm{\d}(w^{-1})\lrcorner X^r(\bm{\lambda})=\overline{\sigma}_-(w)(\d^0\lrcorner \sigma_-(w) X^r(\bm{\lambda})),
\end{equation*}
and in particular we have
$$\bm{\d}(w^{-1})\lrcorner \sigma_+(\mathbf{t}_r)X^r(0)=\overline{\sigma}_-(w)(\d^0\lrcorner \sigma_-(w)\sigma_+(\mathbf{t}_r) X^r(0)).$$
}

\noindent As defined in the introduction, cf. equations \eqref{Erw} and \eqref{Erwtr}, the following two polynomials are required for the next Lemma.

\begin{equation*}
E_r(w)= 1-e_1w+\cdots+(-1)^re_rw^r\in B_r[w]
\end{equation*}
and
\begin{equation*}
E_r\left(\bft_r;\displaystyle{1\over w}\right)=1-e_1(\bft_r)w^{-1}+\cdots+(-1)^re_r(\bft_r)w^{-r}=\prod_{j=1}^r(1-t_jw^{-1})\in \QQ[\bft_r,w^{-1}]
\end{equation*}

\claim{\bf Lemma. }\label{polinomiog} {\em The following identity holds:
 $$\overline{\sigma}_+(w)\overline{\sigma}_{r-1}X^{r-1}(0)\wedge \overline{\sigma}_+(\mathbf{t_r})X^0=\left(E_r(w)+(-1)^{r+1}e_rw^rE_r\left(\bft_r,\displaystyle{1\over w}\right)\right)X^{r}(0).$$
 
}
\medskip

\begin{proof}
We first note that $\overline{\sigma}_{r-1}{X}^{r-1}(0)={X}^{r-1}(1^{r-1})$, and by virtue of \eqref{const_partition} of Lemma \ref{lemconstpart}, where $\ell=1$, we obtain
\begin{equation*}
    \begin{split}
        \overline{\sigma}_+(w)\overline{\sigma}_{r-1}X^{r-1}(0)&=X^{r-1}(1^{r-1})+\sum_{i=1}^{r-1} (-1)^iw^i\,X^{r-1}(2^i1^{r-1-i}),
    \end{split}
\end{equation*}
where for $k\geq \ell$, $(k^{i}\ell^{j})$ stands for the partition having $i$ parts equal to $k$ and $j$ parts equal to $\ell$.
The Proposition \ref{name_4} assures that $\overline{\sigma}_+(\mathbf{t_r})X^0=\sum_{i=0}^{r} (-1)^ie_i(\mathbf{t}_r)X^i$,
and then
\begin{equation*}
    \begin{split}
        X^{r-1}(1^{r-1})\wedge \overline{\sigma}_+(\mathbf{t_r})X^0&=X^{r-1}(1^{r-1})\wedge X^0+(-1)^re_r(\mathbf{t}_r) X^{r-1}(1^{r-1})\wedge X^r\cr\cr
        &=X^r(0)+(-1)^r(-1)^{r-1}e_r(\mathbf{t}_r)e_r X^r(0)\cr\cr
        &=(e_0+(-1)^re_r(-1)^{r-1}e_r(\mathbf{t}_r)) X^r(0).\cr\cr
    \end{split}
\end{equation*}
Now, for each $1\leq i\leq r-1$ we can write
\begin{equation*}
    \begin{split}
        X^{r-1}(2^i1^{r-1-i})=X^{r}\wedge\dots \wedge X^{r+1-i}\wedge X^{r-1-i}\wedge X^{r-2-i}\wedge \dots \wedge X^1,
    \end{split}
\end{equation*}
and so
$$X^{r-1}(2^i1^{r-1-i})\w X^0=X^r(1^i0^{r-i}).$$
Therefore,
\begin{equation*}
    \begin{split}
        &(-1)^iw^i\,X^{r-1}(2^i1^{r-1-i})\wedge \overline{\sigma}_+(\mathbf{t_r})X^0\cr\cr
        &=(-1)^iw^i\,X^{r-1}(2^i1^{r-1-i})\wedge (X^0+(-1)^{r-i}e_{r-i}(\mathbf{t}_r)X^{r-i})\cr\cr
        &=(-1)^iw^i\,X^r(1^i0^{r-i})+(-1)^i(-1)^{r-i}(-1)^{r-i-1}w^ie_{r-i}(\mathbf{t}_r) X^r\wedge \dots\wedge X^1\cr\cr
        &=\left((-1)^ie_iw^i+(-1)^re_r\cdot(-1)^{r-i-1}e_{r-i}(\mathbf{t}_r)w^i\right)X^r(0).\cr
    \end{split}
\end{equation*}
Adding up all the summands above, we get the desired result. \qed
\end{proof}

\medskip

The last required preliminary result right before state and prove the main theorem of this section is the following one.

\claim{\bf Proposition.}\label{name_3} (\cite[Proposition ~4.2]{{GSCH}})
{\em For all $u\in \bigwedge^{r} V$ the following equality holds:
\begin{equation*}
    \sigma_+(z)X^0\wedge u=z^r\sigma_+(z)\overline{\sigma}_-(z)(\overline{\sigma}_ru\wedge X^0).
\end{equation*}
}
\medskip

\claim{\bf Theorem. }\label{thm:mnthm1}
{\em For each $r\geq 2$, $$\Ecal_r(z,w^{-1},\bft_r)=$$
\begin{equation*}
\frac{z^{r-1}}{w^{r-1}}\exp\left(\sum_{n\geq 1}\frac{p_n(\mathbf{t}_r)}{n}\left(\frac{1}{w^n}-\frac{1}{z^n}\right)+x_np_n(z,\mathbf{t}_r)\right)\cdot \left(E_r(w)+(-1)^{r+1}e_rw^rE_r\left(\bft_r,\displaystyle{1\over w}\right)\right).
\end{equation*}
}

\medskip
\claim{\bf Remark.}\label{rmk:rmk412} We call $\Ecal(z,w^{-1},\bft_r)$ the {\em structural formal power series} because it determines, and is determined, by the {\em structural constants} of the representation. If $\frak g$ is an $R$ Lie algebra with basis $(\gamma_a)_{a\in A}$ and  $M$ is a $\frak g$-module with $R$-basis $(m_b)_{b\in B}$, the structural constants of the representation are the scalars $R_{ab}^c$ defined by the equality
$$
\gamma_a\star m_b=\sum_{c\in B}R_{ab}^cm_c\in M. 
$$

\smallskip
\noindent
{\bf Proof of Theorem~\ref{thm:mnthm1}.}
It amounts to work out the  definition \ref{def:def41} of $\Ecal(z,w^{-1},\bft_r)$. One has:
\begin{eqnarray*}
\Ecal_r(z,w^{-1},\bft_r)X^r(0)&=&\sum_{\blamb\in\Pcal_r} \bfX(z)\w (\bmd(w^{-1})\lrcorner X^r(\blamb)s_\blamb(\bft_r))  \hspace{2.65cm}    \cr\cr
&=&\bfX(z)\w (\bmd(w^{-1})\lrcorner\sigma_{+}(\bft_t)X^r(0) \hspace{3.7cm} \mbox{(eq. \eqref{invoke})} \cr\cr
&=&\sigma_+(z)X^0\wedge \overline{\sigma}_-(w)\left(\d^0\lrcorner \sigma_-(w)\sigma_+(\mathbf{t}_r)X^r(0)\right) \ \ \ \ \ \mbox{(Lemma \ref{name})}
\end{eqnarray*}
Now, Proposition \ref{name_1}, equation \eqref{eq:sig-sig+(t_r)}, assures that the expression
$$\sigma_+(z)X^0 \w \ovsig_-(w)\left(\d^0 \lrcorner \sigma_-(w) \sigma_+(\bft_r)X ^r(0)\right)$$
can be written as
$$\exp\left(\sum_{n\geq 1}{1 \over n}p _n\left({\bft_r \over w}\right)\right)\sigma_+(z)X^0 \w \ovsig_-(w)\left(\d^0 \lrcorner  \sigma_+(\bft_r) X ^r(0)\right),$$
Now by Lemma \ref{name_2} this last expression becomes
\begin{equation}\label{maineq1}
    \begin{split}
        &\exp\left(\sum_{n\geq 1}{1 \over n}p _n\left({\bft_r \over w}\right)\right)\sigma_+(z)X^0 \w \ovsig_-(w)\sigma_+(\bft_r)(\d^0 \lrcorner X ^r(0))\\
        &=(-1)^{r-1}\exp\left(\sum_{n\geq 1}{1 \over n}p _n\left({\bft_r \over w}\right)\right)\sigma_+(z)X^0 \w \ovsig_-(w)\sigma_+(\bft_r)X^{r-1}(1^{r-1}).
    \end{split}
\end{equation}
One can see that $\sigma_+(\bft_r)$ commutes with $\ovsig_{-}(w)$ when applied to any $u\in\bigwedge V$ such that $X^0\w u\neq 0\in \bw^{r+1}V$. We have to
avoid $X^0$ just because it is fixed by $\overline\sigma_{-}(w)$.
So in particular,
$\overline{\sigma}_{-}(w)\sigma_+(\bft_r)X^{r-1}(1^{r-1})=\sigma_+(\bft_r)\overline{\sigma}_{-}(w)X^{r-1}(1^{r-1})$.
It is also known that for any $r\geq 1$ one has $\sigma_{+}(\bft_r)\ovsig_{-}(w)X^{r}(1^{r})=\displaystyle{(-1)^{r}\over w^{r}}\ovsig_+(w)\,\sigma_+(\bft_r)X^{r}(0)$.
Hence
\begin{eqnarray*}
\ovsig_-(w) \sigma_+(\bft_r)X^{r-1}(1^{r-1})&=&\displaystyle{(-1)^{r-1} \over w^{r-1}}\ovsig_+(w)\,\sigma_+(\bft_r)X^{r-1}(0).
\end{eqnarray*}

\medskip

%
\noindent The above identity, together Proposition \ref{name_3}, imply that the right hand side of equation \eqref{maineq1} becomes
$${ z^{r-1} \over w^{r-1}}
		\exp\left(\sum_{n\geq 1}{1 \over n}p _n\left({\bft_r \over w}\right)\right)\sigma_+(z)\ovsig_-(z)\left(\ovsig_+(w)\sigma_+(\bft_r)\ovsig_{r-1} X^{r-1}(0)\w X^0\right).$$
Using the fact that $\ovsig_+(w)$ and $\sigma_+(\bft_r)$ commute and integrating by parts \eqref{eq:intp1}
, the previous expression can be written as
\begin{equation}
    \begin{split}\label{maineq2}
        {z^{r-1} \over w^{r-1}}\exp\left(\sum_{n\geq 1}{1 \over n}p _n\left({\bft_r \over w}\right)\right)\sigma_+(z)\ovsig_-(z)\sigma_+(\bft_r)\left(\ovsig_+(w)\ovsig_{r-1}X^{r-1}(0)\w \ovsig_+(\bft_r) X^0\right).
    \end{split}
\end{equation}
To make notation more compact we set
 $$\mathfrak{E}(w,\bft_r):=\left(E_r(w)+(-1)^{r+1}e_rw^rE_r\left(\bft_r,\displaystyle{1\over w}\right)\right).$$
By virtue of Lemma \ref{polinomiog} and also by equation \eqref{eq:sig-sig+(t_r)2} of Proposition \ref{name_1}, the expression in above equation \eqref{maineq2} is equals to
\begin{equation}\label{maineq3}
    \begin{split}
        &{z^{r-1} \over w^{r-1}}\exp\left(\sum_{n\geq 1}{1 \over n}p _n\left({\bft_r \over w}\right)\right)\sigma_+(z)\ovsig_-(z)\sigma_+(\bft_r)\mathfrak{E}(w,\mathbf{t}_r)X^r(0)\\
        &={z^{r-1} \over w^{r-1}}
		\exp\left(\sum_{n\geq 1}{1 \over n}p _n\left({\bft_r \over w}\right)\right)\exp\left(-\sum_{n\geq 1}{1 \over n}p _n\left({\bft_r \over z}\right)\right)
		\sigma_+(z)\sigma_+(\bft_r)\mathfrak{E}(w,\mathbf{t}_r) X^r(0)\\
		&={z^{r-1} \over w^{r-1}}\exp\left(\sum_{n\geq 1}{1 \over n}p _n(\bft_r)\left({1 \over w^n}-{1 \over z^n}\right)\right)\sigma_+(z)\sigma_+(\bft_r)\mathfrak{E}(w,\mathbf{t}_r)X^r(0).
    \end{split}
\end{equation}
Since $X^r(0)$ is eigenvector of $\sigma_+(z)$ with eigenvalue $\displaystyle\frac{1}{E_r(z)}$,  we introduce new formal variables $(x_n)_{n\geq 1}$ through the equality
$$\exp\left(\sum_{n\geq 1}x_nz^n\right)=\frac{1}{E_r(z)}.$$
Hence, substituting 
$$\frac{1}{E_r(z)}\frac{1}{E_r(t_1)}\dots \frac{1}{E_r(t_r)}=\exp\left(\sum_{n\geq 1}x_np_n(z,\bft_r)\right),$$ where $p_n(z,\bft_r)=z^n+p_n(\bft_r)$, 
in last expression of equation \eqref{maineq3} gives 
\begin{equation*}
    {z^{r-1} \over w^{r-1}}\exp\left(\sum_{n\geq 1}{1 \over n}p _n(\bft_r)\left({1 \over w^n}-{1 \over z^n}\right)+x_np_n(z,\bft_r)\right) \mathfrak{E}(w,\bft_r)X^r(0),
\end{equation*}
that concludes the proof of the theorem. \qed

\section{DJKM Representation}\label{sec:sec5}


\claim{} In this section we shall work on the vector space $\Vcal:=\QQ[X^{-1},X]$. We use the following notation:

$$
\wX^0=X^0 \w X^{-1} \w \cdots =X^0 \w \wX^{-1}=X^0 \w X^{-1} \w \wX^{-2}= \cdots.
$$
and if $\blamb=(\lambda_1 , \lambda_2 , \ldots , \lambda_r )\in\Pcal_r$, we set:
$$
\wX^{\blamb}= X^{\lambda_1} \w X^{-1+\lambda_2} \w \cdots \w 
X^{-r+1+\lambda_r} \w \wX^{-r}.
$$
\bclm{\bf Definition.}
{\em The Fermionic Fock space of charge $0$ associated to $\Vcal$ is:}
$$
\Fcal:=\Fcal(\Vcal):=\bigoplus_{\blamb\in\Pcal_r}\QQ\cdot \wX^\blamb.
$$ 
\eclm

\bclm{\bf Definition.} \label{def:extfer}{\em The Schubert derivations extend to $\Fcal$ as follows (see \cite{SDIWP} for details). First one declare that
\begin{eqnarray*}
\sigma_+(z)\wX^m&=&\sigma_+(z)X^m\w \wX^{m-1},\cr\cr
\ovsig_+(z)\wX^m&=&\wX^m-z\wX^{m+1}+z^2\wX^{m+(1^2)}-z^3\wX^{m+(1^3)}+\cdots\,,\cr\cr
\sigma_-(z)\wX^m&=&\wX^m\qquad \mathrm{and}\qquad\ovsig_-(z)\wX^m=\wX^m\label{eq5:extsd}
\end{eqnarray*}
Then one sets:
\begin{eqnarray*}
	\sigma_\pm(z)\wX^{\blamb}&=&\sigma_\pm(z)\left( X^{\lambda_1} \w X^{-1+\lambda_2} \w \cdots \w 
	X^{-r+1+\lambda_r}\right)\w \sigma_\pm(z)\wX^{-r}\cr\cr
\ovsig_\pm(z)\wX^{\blamb}&=&\ovsig_\pm(z)\left( X^{\lambda_1} \w X^{-1+\lambda_2} \w \cdots \w 
    X^{-r+1+\lambda_r}\right)\w \ovsig_{\pm}(z)\wX^{-r}.
\end{eqnarray*}
 }
 \eclm
 
 \medskip

 \bclm{\bf Lemma.}\label{lamma1} {\em 
	The following equation holds:
$$
\sigma_+(t_1)X^0\w \cdots \sigma_+(t_r)X^0 \w\wX^{-r}=\Delta_0(\bft_r)\sigma_+(\bft_r)\wX^0.
$$	

}
\eclm   
\proof Starting on the left hand side, we have

\medskip

\noindent $\sigma_+(t_1)X^0\w \cdots \w \sigma_+(t_r)X^0 \w \wX^{-r} $
\begin{eqnarray*}
	&=&\left(\sum_{i_1,\ldots , i_r \geq 0} X^{i_1}\w \cdots \w X^{i_r} \cdot t_1^{i_1}\cdots  t_r^{i_r}
	\right)\w \wX^{-r}\cr\cr
	&=&\left(\displaystyle{\sum_{(\lambda_1\ldots , \lambda_r)\in\Pcal_r}} X^{\lambda_1}\w \cdots \w X^{-r+1+\lambda_r} s_{\blamb}(\bft_r) \Delta_0(\bft_r)
	\right)\w \wX^{-r}\cr\cr
	&=&\Delta_0(\bft_r)\cdot \sigma_+(t_1,\ldots ,t_r)\big(X^0 \w X^{-1} \w \cdots \w X^{-r+1})\w\wX^{-r}\cr\cr
	&=&\Delta_0(\bft_r)\sigma_+(\bft_r)\big(X^0 \w X^{-1} \w \cdots \w X^{-r+1}\w\wX^{-r}\cr\cr
	&=&\Delta_0(\bft_r)\sigma_+(\bft_r)\wX^0\hskip 150pt
\end{eqnarray*} as desired. \qed

\bclm{\bf Lemma.} \label{lemma2}{\em We have:
	$$
	\displaystyle{\sum_{\bmu\in\Pcal_k}}\wX^{\blamb} s_{\blamb}(\bft_r) \w \wX^{-r}=\sigma_+(\bft_r)\wX^0.
	$$	
}
\proof 
By \cite[Proposition 5.12]{SDIWP} the product $\sigma_+(t_1)\cdots\sigma_+(t_r)$ of $r$ Schubert derivations acts only on the first $r$ exterior factors of $\wX^0$, namely
$$
\sigma_+(\bft_r)\wX^0=\sigma_+(\bft_r)\big(X^0\w\cdots\w X^{-r+1}\big)\w \wX^{-r}.\hskip 100pt \qed
$$

\claim{\bf Lemma.}\label{lamma3} {\em The following commutation rule holds:
$$
\sigma_-(z)\sigma_+(t)\wX^0 =\exp\left(\sum_{n  \geq 1} {1 \over n}{t^n \over z^n}
\right)   \sigma_+(t)\sigma_-(z)\wX^0=\exp\left(\sum_{n  \geq 1} {1 \over n}{t^n \over z^n}
\right)   \sigma_+(t)\wX^0.
$$
}
\eclm
\proof It amounts to straightforward mechanical manipulation which we report below with details just for sake of completeness. Using Definition \ref{def:extfer} for the extension of Schubert derivations to the Fermionic Fock space
\begin{eqnarray*}
\sigma_-(z)\left[\sigma_+(t)\wX^0\right]&=&\sigma_-(z)\left(
\sigma_+(t)X^0 \w \wX^{-1} \right)\cr\cr
&=&\sigma_-(z)\sigma_+(t)X^0 \w \wX^{-1}=\sigma_-(z)\left(\sum_{i\geq 0} X^i t^i\right)\w \wX^{-1}\cr\cr
&=&\left[
\left(X^0+ {X^{-1} \over z}+{X^{-2} \over z^2}+\cdots\right)+
\left(X^1+ {X^{0} \over z}+{X^{-1} \over z^2}+\cdots\right)t\right.\cr
&& \left.
+\left(X^2+ {X^{1} \over z}+{X^{0} \over z^2}+\cdots\right)t^2 
+\cdots \right]\w X^{-1} \w X^{-2}\w \cdots\cr\cr
&=&\left[
\left(X^0 +X^1 t+ X^2 t^2 + \cdots\right)+
{t \over z}\left(X^0 +X^1 t+ X^2 t^2 + \cdots\right)\right.\cr
&& \left.
+{t^2 \over z^2}\left(X^0 +X^1 t+ X^2 t^2 + \cdots\right)+\cdots
\right]\w \wX^{-1}\cr\cr
&=& \sigma_+(z)X^0\left(1+{t \over z}+{t^2 \over z^2}+\cdots\right)\w \wX^{-1}\cr\cr
&=&\exp\left(\sum_{n  \geq 1} {1 \over n}{t^n \over z^n}
\right)  \sigma_+(z)X^0 \w \wX^{-1}\cr\cr
&=&\exp\left(\sum_{n  \geq 1} {1 \over n}{t^n \over z^n}
\right)  \sigma_+(z)\wX^{0}\cr
\end{eqnarray*}
\qed

\bclm{\bf Lemma.}\label{lamma4} {\em  The following commutation rule holds
	$$
	\ovsig_-(z)\sigma_+(t)\wX^0 =\exp\left(-\sum_{n  \geq 1} {1 \over n}{t^n \over z^n}
	\right)   \sigma_+(t)\ovsig_-(z)\wX^0=\exp\left(-\sum_{n  \geq 1} {1 \over n}{t^n \over z^n}
	\right)   \sigma_+(t)\wX^0
	$$
}
\eclm
\proof First one writes:
$$
\sigma_+(t)\wX^0=\ovsig_-(z)\sigma_-(z)\sigma_+(t)\wX^0.
$$
Now we use Lemma \ref{lamma3}, to commute $\sigma_-(z)$ and $\sigma_+(t)$, so obtaining:

\begin{eqnarray*}
\sigma_+(t)\wX^0&=& \exp\left(\sum_{n  \geq 1} {1 \over n}{t^n \over z^n}
\right) \ovsig_-(z)\sigma_+(z)\wX^0,
\end{eqnarray*}
from which 
$$
\hskip50pt \ovsig_-(z)\sigma_+(t)\wX^0=\exp\left(-\sum_{n  \geq 1} {1 \over n}{t^n \over z^n}
\right)   \sigma_+(z)\wX^0.\hskip 160pt \qed
$$

Recall the notation for Vandermonde determinant in the variables $\bft_r$ as in \eqref{eq:vande}.
\bclm{\bf Corollary.}\label{Corollary5} {\em 
For all $r \geq 1$,
$$
\sigma_-(z)\sigma_+(\bft_r)\wX^0=
\exp\left(\sum_{n  \geq 1} {1 \over n}{{p_n(\bft_r)} \over z^n}
\right)   \sigma_+(\bft_r)\wX^0.
$$	
}
\eclm
\proof One has:
\begin{eqnarray*}
\Delta_0(\bft_r)\sigma_-(z)\sigma_+(\bft_r)\wX^0&=&
\sigma_-(z)\left(\sigma_+(t_1)X^0 \w \cdots \w \sigma_+(t_r)X^0\right)\w \wX^{-r}\cr\cr
&=&\Delta_0(\bft_r)\sigma_-(z)\sigma_+(t_1)X^0  \w \cdots 
\sigma_-(z)\sigma_+(t_r)X^0  \w \wX^{-r}\cr\cr
&=&\Delta_0(\bft_r) 
\exp\left(\sum_{n  \geq 1} {1 \over n}{t_1 ^n \over z^n}\right)\cdots 
\exp\left(\sum_{n  \geq 1} {1 \over n}{t_r ^n \over z^n}
\right)
\sigma_+(\bft_r)\wX^0\cr\cr
&=&\Delta_0(\bft_r) \exp\left(\sum_{n  \geq 1} {1 \over n}{{p_n(\bft_r)} \over z^n}
  \right)   \sigma_+(\bft_r)\wX^0
\end{eqnarray*}
The proof ends by dividing the very first and the very last member by $\Delta_0(\bft_r)$.\qed
\bclm{\bf Corollary.}\label{Corollary5.2} {\em 
	For all $r \geq 1$,
	$$
	\ovsig_-(z)\sigma_+(\bft_r)\wX^0=
	\exp\left(-\sum_{n  \geq 1} {1 \over n}{{p_n(\bft_r)} \over z^n}
	\right)   \sigma_+(\bft_r)\wX^0.
	$$	
}
\eclm
\proof It is an  obvious consequence of Corollary \ref{Corollary5}, playing with the inverse of the Schubert derivations.\qed
 
The Boson--Fermion correspondence implies that $B=\QQ[x_1 , x_2 , \dots ]$ is naturally isomorphic to $\Fcal_0$ via the $\QQ$-linear extension of the sets map 
\begin{equation*}
S_{\blamb}(\bfx)\overset{\cong}{\longmapsto} \wX^{\blamb}.\label{eq1:fbfc}
\end{equation*}
\claim{} Gatto \& Salehyan in \cite{SDIWP} show that 
$$
\left(
\exp \left(\sum_{n\geq 1}(z^n -w^n)\right)S_{\blamb}(\bfx)
\right)\wX^0=\sigma_+(z)\ovsig_+(w)\wX^\blamb,
$$
and that
$$
\left(
\exp \left(-\sum_{n\geq 1}{1\over n}\left({1 \over z^n} -{1 \over w^n}\right){\d\over \d x_n}\right)S_{\blamb}(\bfx)
\right)\wX^0=\ovsig_-(z)\sigma_-(w)\wX^{\blamb}.
$$
As in the case for finite $r$,  $gl(\Vcal)=\Vcal \otimes \Vcal^*$ acts on $B$ as follows:
$$
\left[
\left(
X^i \otimes \d ^j \right)S_{\blamb}(\bfx)
\right]\wX^0=X^i \w \d ^j \lrcorner \wX^{\blamb}.
$$
Putting as usual $$\Ecal(z,w)=\sum_{i.j \in \ZZ}\left(X^i \w \d ^j \lrcorner\right)z^i w^{-j},$$
The DJKM result \cite{DJKM01} (see also \cite[Theorem 5.1]{KacRaRoz}) says that:
\be
\Ecal(z,w)={1\over 1-\displaystyle{w\over z}} \cdot \Gamma (z,w),  \label{formula1}
\ee
where
\begin{equation*}
\Gamma (z,w)=\exp \left(\sum_{n\geq 1}(z^n -w^n)\right)
\exp \left(-\sum_{n\geq 1}{1\over n}\left({1 \over z^n} -{1 \over w^n}\right){\d\over \d x_n}\right).\label{formula2}
\end{equation*}
\bclm{\bf Theorem.}\label{thm:mnthm2} {\em
Let 
$$\Ecal(z,w^{-1}, \bft_r):=\displaystyle{\sum_{\blamb \in\Pcal_r}} \Gamma (z,w) S_{\blamb}(\bfx)\bfs_{\blamb}(\bft_r)\in B \llb z,w, \bft_r,z^{-1}, w^{-1}\rrb 
$$ 
be the {\em structural formal power series} of the $B$-representation of $gl(\Vcal)$ (Cf. Remark~\ref{rmk:rmk412}).
Then 
$$
\Ecal(z,w^{-1}, \bft_r)=
\exp \left(\sum_{n\geq 1}{1\over n}\left({w^n \over z^n} -{p_n (\bft_r) \over z^n}+{p_n (\bft_r) \over w^n}\right)
+x_n \left(z^n -w^n+p_n (\bft_r)\right)
\right).
$$
}
\eclm
\proof We first use the exponential expression of the geometric series:
$$
\left(1-\displaystyle{w\over z}\right)^{-1}=
\exp\left(\sum_{n  \geq 1} {1 \over n}{w ^n \over z^n}\right),
$$
to put in \eqref{formula1}.
Now we observe that:
\begin{eqnarray*}
\left(\displaystyle{\sum_{\blamb \in\Pcal_r}} \Gamma (z,w) S_{\blamb}(\bfx)\bfs_{\blamb}(\bft_r)\right)\wX^0&=&
\sigma_+(z)\ovsig_+(w)\ovsig_-(z)\sigma_-(w)
\displaystyle{\sum_{\blamb \in\Pcal_r}}\wX^{\blamb} \bfs_{\blamb}(\bft_r)
\cr\cr
&=&\sigma_+(z)\ovsig_+(w)\ovsig_-(z)\sigma_-(w)\sigma_+(\bft_r)\wX^0.
\end{eqnarray*}
Applying Lemma \ref{lamma3} and Lemma \ref{lamma4}:
\begin{eqnarray}
&=&	\exp\left(\sum_{n  \geq 1} {1 \over n}{{p_n(\bft_r)} \over w^n}\right) 
\sigma_+(z)\ovsig_+(w)\ovsig_-(z)\sigma_+(\bft_r)\wX^0	\cr\cr
&=& \exp\left(\sum_{n  \geq 1} {1 \over n}{{p_n(\bft_r)} \over w^n}\right) 
\exp\left(-\sum_{n  \geq 1} {1 \over n}{{p_n(\bft_r)} \over z^n}\right) 
\sigma_+(z)\ovsig_+(w)\sigma_+(\bft_r)\wX^0\label{simpl1}.
\end{eqnarray}	
By recalling that $\wX^0$ is eigenvalue of $\sigma_+(z),\ovsig_+(w)$ and $\sigma_+(\bft_r)$, we obtain
\begin{equation*}\label{simpl2}
\sigma_+(z)\ovsig_+(w)\sigma_+(\bft_r)\wX^0= \exp \left(\sum _{n\geq 1}x_n z^n\right)\exp \left(-\sum _{n\geq 1}x_n w^n\right) \exp \left(\sum _{n\geq 1}x_n p_n (\bft_r)\right)\wX^0.
\end{equation*}
Substituting the previous equation on \eqref{simpl1} and simplifying, we may conclude
\begin{eqnarray*}
\Ecal(z,w^{-1}, \bft_r)&=&\displaystyle{\sum_{\blamb \in\Pcal_r}} \Gamma (z,w) S_{\blamb}(\bfx)\bfs_{\blamb}(\bft_r)\cr\cr
&=&\exp \left(\sum_{n\geq 1}{1\over n}p_n (\bft_r)\left(\displaystyle{1 \over z^n}-\displaystyle{1 \over w^n}\right)
+x_n \left(z^n -w^n+p_n (\bft_r)\right)
\right),
\end{eqnarray*}	
as claimed. \qed

\bibliographystyle{amsplain}
\providecommand{\bysame}{\leavevmode\hbox to3em{\hrulefill}\thinspace}
\providecommand{\MR}{\relax\ifhmode\unskip\space\fi MR }
\providecommand{\MRhref}[2]{%
  \href{http://www.ams.org/mathscinet-getitem?mr=#1}{#2}
}
\providecommand{\href}[2]{#2}
\bibliographystyle{amsplain}


\medskip
\medskip

\parbox[t]{3in}{{\rm Ommolbanin~Behzad}\\
{\tt \href{mailto:O.behzad@mcs.ui.ac.ir}{O.behzad@mcs.ui.ac.ir}}\\
{\it Department of Pure Mathematics}\\
{\it Faculty of Mathematics and Statistics}\\
{\it University of Isfahan}\\
{\it P.O.Box  81746-73441,  Isfahan}\\
{\it  IRAN}} \hspace{1.5cm}
\parbox[t]{3in}{{\rm Andr\'e Contiero}\\
{\tt \href{mailto:contiero@ufmg.br}{contiero@ufmg.br}}\\
{\it Universidade Federal de Minas Gerais}\\
{\it Belo Horizonte, MG}\\
{\it BRAZIL}}

\bigskip

\vspace{6 pt}

\parbox[t]{3in}{{\rm David Martins}\\
{\tt \href{mailto:davidtm@ufmg.br}{davidtm@ufmg.br}}\\
{\it Universidade Federal de Minas Gerais}\\
{\it Belo Horizonte, MG}\\
{\it BRAZIL}} \hspace{1.5cm} 

\end{document}